\documentclass{elsart}

\usepackage{amsmath}
\usepackage{amssymb}
\journal{Journal of Number Theory}

\newcommand{\Primes}{{\mathbb P}}

\newcommand{\floor}[1]{{\lfloor #1 \rfloor}}
\newcommand{\h}{h}

\def\cocoa{{\hbox{\rm C\kern-.13em o\kern-.07em C\kern-.13em o\kern-.15em A}}}

\begin{document}
\begin{frontmatter}

\title{On the Number of Distinct Multinomial Coefficients}

\author{George E. Andrews\thanksref{GEA}\corauthref{cor}}
\address{Mathematics Department, 410 McAllister Building,
 The Pennsylvania State University, University Park, PA 16802,
 USA.}
\ead{andrews@math.psu.edu} \corauth[cor]{Corresponding author.}
\thanks[GEA]{Partially supported by
National Science Foundation Grant DMS-0200047.}

\author{Arnold  Knopfmacher\thanksref{AK}}
\address{The John Knopfmacher Centre for Applicable Analysis and Number Theory,
University of the Witwatersrand, Johannesburg, South Africa.}
\ead{arnoldk@cam.wits.ac.za}
\thanks[AK]{Partially  supported by The John Knopfmacher Center for
Applicable Analysis and  Number  Theory  of  the  University  of
the Witwatersrand.}

\author{Burkhard  Zimmermann\thanksref{BZ}}
\address{Research Institute for Symbolic Computation, Johannes Kepler Universit\"{a}t Linz, A-4040 Linz,
Austria.} \ead{Zimmermann@risc.uni-linz.ac.at}
\thanks[BZ]{Supported  by SFB grant F1301 of the
Austrian FWF.}
\begin{abstract}
We  study  $M(n)$, the number of distinct values
taken  by  multinomial  coefficients with upper entry $n$,
and some closely related sequences.
We show that  both  $p_{\Primes}(n)/M(n)$ and $M(n)/p(n)$ tend to zero as
$n$  goes  to  infinity,  where $p_{\Primes}(n)$ is the number of
partitions  of  $n$  into primes and $p(n)$ is the total number of
partitions of $n$. To use methods from commutative algebra,
we encode partitions and multinomial coefficients as monomials.
\end{abstract}


\end{frontmatter}

\section{Introduction}
The classical multinomial expansion is given by
\begin{equation} \label{eq:expansion}
(x_1 + x_2 + \cdots + x_k)^n =
\sum \binom{n}{i_1,i_2,\dots,i_k} x_1^{i_1}x_2^{i_2} \cdots x_k^{i_k}\,,
\end{equation}
where the sum
runs over all $(i_1,i_2,\dots,i_k)$ such that
$i_1 + i_2 + \cdots + i_k = n$ and $i_1,i_2,\dots,i_k \geq 0$.
Multinomial coefficients are defined by
\begin{equation} \binom{n}{i_1,i_2,\dots,i_k}
:= \frac{n!}{i_1! i_2! \cdots i_k!}. \end{equation}

It  is  natural to  ask  about $M_k(n)$,
the number of different values of
\begin{equation} \binom{n}{i_1,i_2,\dots,i_k}\end{equation}
where  $i_1+i_2+\dots+i_k=n$. Obviously if the $i_1,i_2,\dots,i_k$
are  merely  permuted, then the value of $\binom{n}{i_1,i_2,\dots,
i_k}$ is  unchanged.  However  identical values do not necessarily
arise  only  by  permuting  the  $i_1,i_2,\dots,i_k$. For example,
\begin{equation} \label{eq:3224111}
\binom{7}{3,2,2} = \binom{7}{4,1,1,1}\end{equation}
and
\begin{equation}
\binom{236}{64,\,55,\,55,\,52,\,7,\,3}
= \binom{236}{62,\,56,\,54,\,51,\,13}.
\end{equation}

We  note that if $k \geq n$, then $M_k(n) = M_n(n)$, and we define
$M(n)  :=  M_n(n)$  to be the total number of distinct multinomial
coefficients with upper entry~$n$.

Since   permuting   its  lower  indices  leaves  the  value  of  a
multinomial  coefficient  unchanged  it  is immediately clear that
\begin{equation} \label{ineq:Mkpk}
M_k(n) \leq p_k(n) \end{equation}
and
\begin{equation}
M(n)  \leq  p(n)\,, \end{equation} where $p_k(n)$ is the number of
partitions  of $n$ into at most $k$ parts, and $p(n)$ is the total
number  of  partitions  of  $n$  respectively.  Observing that the
binomial coefficients $\binom{n}{k,\,n-k}$ are strictly increasing
for  $0  \leq  k  \leq  \frac{n}{2}$,  we  deduce  that,  in fact,
\begin{equation} \label{eq:M2p2}
M_2(n) = p_2(n).
\end{equation}
However  the inequality (\ref{ineq:Mkpk}) seems to be stronger for
large $k$. Indeed (Theorem~\ref{thm:Mp0}),
\begin{equation}
\lim_{n\rightarrow \infty} \frac{M(n)}{p(n)} = 0.
\end{equation}
Bounding $M(n)$ from below we will prove
(Theorem~\ref{lower-bound-thm}) that
\begin{equation}M(n) \geq p_{\Primes}(n)\end{equation}
where  $p_{\Primes}(n)$
is  the number of partitions of $n$ into
parts belonging to the set of primes
$\Primes$. Indeed
(Theorem~\ref{thm:pM0}),
\begin{equation}\lim_{n\rightarrow \infty}
\frac{p_{\Primes}(n)}{M(n)} = 0.\end{equation}

It  is  natural to generalize the problem from $M(n)$ to $M_S(n)$,
the  number of different multinomial coefficients with upper entry
$n$  whose  lower  entries  belong  to  a  given set $S$ of
natural  numbers.
Let
\begin{equation}\label{eq:Mgen}
  \mathcal{M}_S(q) := \sum_{n} M_S(n) q^n
\end{equation}
and
\begin{equation}\label{eq:pgen}
  \mathcal{P}_S(q) := \sum_{n} p_S(n) q^n
\end{equation}
where $p_S(n)$ is the number of partitions of $n$ into elements
from $S$.
Define $[s] := \{1,2,\dots,s\}$. Results of numerical calculations
such as
\begin{equation}\label{eq:45}
  \mathcal{M}_{[4]}(q) \,/\, \mathcal{P}_{[4]}(q)
= 1 - q^7 + O(q^{100})
\end{equation}
and
\begin{equation}\label{eq:67}
  \mathcal{M}_{[7]}(q) \,/\, \mathcal{P}_{[7]}(q)
= 1 - q^7 - q^8 -q^{10} + q^{12} + q^{13} + O(q^{100})
\end{equation}
suggest that
$\mathcal{M}_{S}(q) \,/\, \mathcal{P}_{S}(q)$
is a polynomial for any finite $S$.
This is indeed true (Theorem \ref{thm:form}) and leads
to an algorithm
for computing a closed form for the sequence $M_S(n)$
for a given finite set $S$ (Section~\ref{sec:finite}).

Partitions and multinomial coefficients can be written as monomials in
a natural way:
For instance, the monomial $q_4 q_1^3$ represents the partition $4+1+1+1$, and
$x_7 x_5 x_3 x_2$ represents the multinomial coefficient $\binom{7}{4,1,1,1}$
whose factorization into primes is $7 \cdot 5 \cdot 3 \cdot 2$.
This encoding serves as a link between our counting problem
and Hilbert functions (Section \ref{sec:setting}). Sections \ref{sec:finite},
\ref{sec:uppers} and \ref{sec:lowers} are based on that link.

We  call  a  pair  of  partitions  of  $n$  that  yield  the  same
multinomial   coefficient   but  have  no  common  parts  an  {\em
irreducible  pair}.
For  example,  the  partitions  $4+1+1+1$ and
$3+2+2$     form    an    irreducible    pair    according    to
Equation~(\ref{eq:3224111}).
In  Section \ref{sec:in}, we study
$i(n)$ the total number of irreducible pairs of partitions of $n$,
and  we  prove (Theorem \ref{thm:n56}) that $i(n) > \frac{n}{56} -
1$.

\section{A Lower Bound for $M(n)$} \label{sec:lowerbound}
We relate $M(n)$ to $p_{\Primes}(n)$ whose asymptotics is
known   by   a   theorem   of Kerawala \cite{Ker69}:
\begin{equation} \log p_{\Primes}(n) \sim
\frac{2 \pi}{\sqrt{3}} \sqrt{\frac{n}{\log  n}}.
\end{equation}
\begin{thm} \label{lower-bound-thm}
$M(n) \geq p_{\Primes}(n)$.\end{thm} Theorem~\ref{lower-bound-thm}
is implied by the following lemma:
\begin{lem}\label{NeverEquivalent}
Any two distinct partitions of the same natural number $n$
into  primes  yield  different  multinomial
coefficients.
\end{lem}
\begin{pf}[Proof of Lemma \ref{NeverEquivalent}]
It suffices to show that if
\begin{equation}
p_1!  p_2!  \cdots  p_r!  =  q_1!  q_2! \cdots q_s! \end{equation}
where $p_1 \leq p_2 \leq \cdots \leq p_r$ and $q_1 \leq q_2 \leq \cdots
\leq q_s$  are  all  primes, then $r = s$ and $p_i = q_i$ for $i = 1,
\ldots, s$. We proceed by mathematical induction on $r$.

If  $r  = 1$, then $q_s$ must equal $p_1$  because  if $q_s < p_1$
then $p_1$ divides the left side of the above equation but not the
right side.  If $q_s > p_1$  then $q_s$ divides the right side but
not  the  left.  Hence  $q_s  =  p_1$,  and dividing both sides by
$p_1!$  we  see  that there can be no other $q_i$. Hence $s=1$ and
$q_1 = p_1$.

Assume  now  that  our  result  holds  up  to  but not including a
particular  $r$.  As  in the case $r=1$, we must have $q_s = p_r$.
Cancel  $p_r!$  from both sides and apply the induction hypothesis
to  conclude  that $s-1 = r-1$ and $p_i = q_i$ for $i = 1, \ldots,
s-1$. Hence the lemma follows by mathematical induction.
\qed\end{pf}
Some  values  of  $p_{\Primes}(n)$  and $M(n)$ are listed on page
\pageref{numbers}.
We  will  refine Theorem~\ref{lower-bound-thm} in
Section    \ref{sec:lowers}.

\section{The Algebraic Setting}\label{sec:setting}

Encoding partitions and multinomial coefficients as monomials
allows us to apply constructive methods from commutative algebra
to the problem of counting multinomial coefficients.
Let us assume that $S \subseteq \mathbb{N}$ throughout the paper.
We will see that $M_S(n)$ finds a natural interpretation as the Hilbert
function of a certain graded ring (Lemma~\ref{lem:hilbert}).
In the case of finite $S$, it can be computed by
the method of Gr\"{o}bner bases \cite{AL94,BB70,BB85,CLO92}.

We represent the partition $\lambda_0 + \lambda_1 + \dots + \lambda_i$
of $n$ by the monomial
$q_{\lambda_0} q_{\lambda_1} \dots q_{\lambda_i}$
whose degree is $n$ if we define the degree of variables
suitably by $\deg q_j := j$.
For convenience, we will use the notions
``partition of $n$'' and ``monomial of degree $n$''
interchangeably.

Let  $k$ be a field of characteristic zero. We abbreviate the ring
$k[q_i : i \in S]$ of polynomials  in  the  variables  $q_i$ for
$i \in S$ over $k$ by $k[S]$. Define  the  degree  of  monomials
by $\deg  q_i  := i$, and let $k[S]_n$ denote the subspace of all
homogeneous polynomials of degree $n$. In other words, $k[S]_n$ is
the $k$-vector space whose basis are the partitions of $n$ into
parts $S$. Note that $k[S]$ is graded by $k[S] = \bigoplus_n
k[S]_n$. For instance,
\begin{equation}
k[\{1,\, 3,\dots \}]_4 \,=\,
k \cdot q_3 q_1 \oplus\, k \cdot {q_1}^4
\end{equation}
corresponding  to  the  partitions
$3+1$ and $1+1+1+1$
of $4$ into odd parts.

The  multinomial  coefficients  with  upper  entry  $n$ into parts
belonging  to  $S$  are  the numbers $n! / \prod_j j!^{a_j}$ where
$\prod_j  q_j^{a_j}$  ranges  over the monomials in $k[S]_n$. Since
the  numerator  $n!$  of  these fractions is fixed, it suffices to
count the set of all denominators:
\begin{equation} \label{eq:prodj} M_S(n) = |\{
\prod_j j!^{a_j} \,:\, \prod_j q_j^{a_j} \in k[S]_n \}|.
\end{equation}
To  count  the  values taken by $\prod_j j!^{a_j}$, we look at
their factorization  into primes. Let $\h(q_j)$ be the
factorization of $j!$ into primes, written as a monomial in $k[x]
:= k[x_p\,:\,p \; \text{prime}]$, multiplied by $q^j$. For
example, $ \h(q_5) = q^5 x_2^3 \, x_3 \, x_5$ corresponding to $5!
= 2^3 \cdot 3 \cdot 5$. An elementary counting argument
\cite{GKP94} shows that the prime $p$ occurs in the factorization
of $j!$ with exponent $\sum_{l=1}^\infty  \floor{j/p^l}$, where
$\floor{x}$ denotes the largest  integer  that  does  not exceed
the real number $x$. Therefore,
\begin{equation}\label{eq:phidef}           \h(q_j)           =
q^j  \prod_{p \; \text{prime}} {x_p}^{\sum_{l=1}^\infty \floor{j/p^l}}.\end{equation}
Since factorization into primes is unique, (\ref{eq:prodj})
can be written as
\begin{equation}\label{eq:MSphi} M_S(n) = |\{
\prod_j \h( q_j )^{a_j} \,:\, \prod_j q_j^{a_j} \in k[S]_n \}|.
\end{equation}
Extending  $\h$ to a $k$-algebra homomorphism $k[S] \rightarrow
k[x,q]$ allows  us to reformulate (\ref{eq:MSphi}) as
\begin{lem}\label{lem:dim}
\begin{equation} \label{eq:dimphi}
M_S(n) = \dim_k \h(k[S]_n). \end{equation}
\end{lem}

{\em Example:} Since there are $10$ partitions of $7$ into parts $1,2,3$
and $4$, the dimension of
\begin{equation}k[\{1,2,3,4\}]_7 =
       k \, q_4 q_3
\oplus k \, q_4 q_2 q_1
\oplus k \, q_4 q_1^3
\oplus k \, q_3 q_2^2
\oplus \dots
\oplus k \, q_1^7\end{equation}
over $k$ is $10$.
However, the dimension of its image
\begin{equation}
\h( k[\{1,2,3,4\}]_7 ) =
       k \, q^7 x_2^4 x_3^2
\oplus k \, q^7 x_2^4 x_3
\oplus k \, q^7 x_2^3 x_3
\oplus \dots
\oplus k \, q^7
\end{equation}
under $\h$ is only $9$ and so $M_{[4]}(7)=9$. The defect is due to
$\h(q_4 q_1^3) = \h(q_3 q_2^2)$ which is nothing but a restatement
of (\ref{eq:3224111}).

To use Lemma~\ref{lem:dim} for effective computation (in the case
of finite $S$), we express $\dim_k \h(k[S]_n)$ as the value (at
$n$) of the Hilbert function of a certain elimination ideal. This
method is taken from \cite{AL94}; the result in our case is
Lemma~\ref{lem:hilbert} below.

First we make the map $\h$ degree-preserving (graded) by defining
$\deg q := 1$ and $\deg x_p := 0$ in the ring $k[x_p \,:\, p \;
\text{prime}][q]$. (This is why we introduced the extra factor of
$q^j$ in the defining equation (\ref{eq:phidef}) of $\h$.) Second,
note that
\begin{equation}
\h(k[S]_n) \cong k[S]_n / (k[S]_n \cap \,\ker \h)\end{equation} as
$k$-vector  spaces, since $\h$ is a $k$-linear map on $k[S]_n$. In
particular, dimensions agree. Therefore,
\begin{equation} \label{eq:dimquot}
M_S(n) = \dim_k k[S]_n / (k[S]_n \cap \,\ker \h).
\end{equation}
Recall  that  the  (projective) Hilbert function $H_R$ of a graded
$k$-algebra  $R = \bigoplus_n R_n$ is defined by $H_R(n) := \dim_k
R_n$.   Thus  (\ref{eq:dimquot})  relates  $M_S$  to  the  Hilbert
function  of  $k[S]  /  \ker \h$:
\begin{equation} \label{eq:hilbertker}
M_S(n) = H_{k[S] / \ker \h}(n).
\end{equation}
By Theorem  2.4.2  of  \cite{AL94}, $\ker \h$ can be computed by
elimination:
\begin{equation} \ker \h = I \cap k[S]\end{equation}
where the  ideal  $I$  of  $k[S][q][x_p \,:\, p \,\text{prime}]$
is defined by
\begin{equation} I :=
\langle q_j - \h(q_j) \,:\, j\in S \rangle.\end{equation}
Summarizing this section, we have proved the following Lemma:
\begin{lem}\label{lem:hilbert}
Let $k[S]$ be graded by $\deg  q_i  := i$. Define a $k$-algebra
homomorphism from $k[S]$  to $k[q,x]$ by
\begin{equation}\label{eq:phidef1}           \h(q_j)           :=
q^j  \prod_p {x_p}^{\sum_{l=1}^\infty
\floor{j/p^l}}.\end{equation} Let  the  ideal  $I$  of $k[S,q,x]$
be defined by
\begin{equation} I :=
\langle q_j - \h(q_j) \,:\, j\in S \rangle\end{equation} and let
\begin{equation} J := I \cap k[S].\end{equation}
Then $M_S$ is the (projective) Hilbert function of the
$k$-algebra $k[S] / J$:
\begin{equation} \label{eq:hilbert}
M_S(n) = H_{k[S] / J}(n).
\end{equation}
\end{lem}
{\em Example:}
If $S=[4]$, then
$I = \langle
q_1 - q, \,
q_2 - q^2 x_2, \,
q_3 - q^3 x_2 x_3, \,
q_4 - q^4 x_2^3 x_3 \rangle$
and $J = \langle q_4 q_1^3 - q_3 q_2^2 \rangle$.
For $M_{[4]}(n)$, see (\ref{eq:M4n}) on page~\pageref{eq:M4n}.

\section{Explicit Answers}\label{sec:finite}

Let $S$ be a given finite set throughout this section.
Lemma~\ref{lem:hilbert} allows to compute a closed form for the
sequence $M_S(n)$ by well-known methods from computational
commutative algebra. For the sake of completeness, let us briefly
review them:

\begin{enumerate}
\item    Fix a    term    order   $\preceq$   on $k[S,q,x]$ that
allows the elimination of  the  variable  $q$  and the variables
$x_p$ in step 2 below. Compute a Gr\"{o}bner basis $F$ for the
(toric) ideal $I = \langle   q_j  -  \h(q_j)  \,:\,  j\in  S
\rangle$ with  respect  to  this  term  order  using Buchberger's
algorithm \cite{BB70,BB85}. \item Let   $G  :=  F  \cap  k[S]$. By
the elimination  property  of Gr\"{o}bner  bases with respect to a
suitable elimination order $\preceq$, the set $G$ is a Gr\"{o}bner
basis for the elimination ideal $ J = I \cap k[S]$. \item Let  $L
:=  I_\preceq  (G)$  be  the set  of  leading  terms of
polynomials  in  $G$. \item Compute $\mathcal{M}_S(q)$ using
\begin{equation}
\mathcal{M}_S(q)
= \mathcal{H}_{k[S] / J}(q)
= \mathcal{H}_{k[S]   /  I_\preceq  (J)}(q)
= \mathcal{H}_{k[S]   /  \langle L \rangle}(q).  \end{equation}
The first equality holds by Lemma~\ref{lem:hilbert}. The second
equality is an identity of Macaulay~\cite{Mac27}.
Since  $G$  is a  Gr\"{o}bner  basis, its initial terms $L$ generate
the initial term ideal of $\langle G \rangle$ with respect to $\preceq$,
which explains the third equation sign.
A naive method
for computing the Hilbert-Poincar\'{e} series of
$k[S] /\langle L \rangle$
is to apply the inclusion-exclusion relation
\begin{equation}\label{eq:HtL}
             \mathcal{H}_{k[S] / \langle \{t\} \cup L \rangle}(q)
=
             \mathcal{H}_{k[S] / \langle L \rangle}(q)
- q^{\deg t} \mathcal{H}_{k[S] / \langle L \rangle : t}(q),
\end{equation}
recursively until the base case
\begin{equation}\label{eq:Hbase}
  \mathcal{H}_{k[S] / \langle \rangle}(q)
= \mathcal{H}_{k[S]}(q) = \frac{1}{\prod_{j \in S} (1 - q_j)}
\end{equation}
is  reached. For better (faster) algorithms, see~\cite{B97}.
\item Extract a closed form expression for
$H_{k[S]  /  \langle  L  \rangle}(n)$ from its
generating function $ \mathcal{H}_{k[S] / \langle L \rangle}(q)$.
(Use partial fraction decomposition and the binomial series).
It is the desired answer $M_S(n)$.
\end{enumerate}
One of the authors computed 1 -- 4 for several finite $S$ using
different computer algebra systems. It turned out that \cocoa
\cite{CocoaSystem} was fastest for that purpose.
\begin{thm}\label{thm:form}
Let $S$ be a finite subset of the positive natural numbers.
Then
\begin{enumerate}
\item $\mathcal{M}_S(q)$ can be written as
\begin{equation}\label{eq:form}\mathcal{M}_S(q)
= \frac{f_S(q)}{\prod_{j\in S} (1-q^j)}\end{equation}
where $f_S(q)$ is a polynomial with  integer  coefficients.
\item
There exists $n_0$ such that
$M_S(n)$ can be written as a quasipolynomial~\cite{St}
for $n \geq n_0$.
Moreover, it suffices to use periods which are divisors of elements of $S$.
\end{enumerate}
\end{thm}
\begin{pf}
Relations (\ref{eq:HtL}) and (\ref{eq:Hbase}) prove the first statement.
The second statement follows from the first easily.
\qed\end{pf}

Let  us  follow the algorithm in the case $S = [4]$, which
is the simplest nontrivial case. We have $ I = \langle q_1 - q,\,
q_2  -  q^2  x_2,\,  q_3  -  q^3  x_2  x_3,\,  q_4 - q^4 x_2^3 x_3
\rangle$. To
eliminate  the  variables  $x_3$,  $x_2$ and $q$ we
choose  the  lexical term order where $x_3 \succ x_2 \succ q \succ
q_4  \succ  q_3  \succ  q_2  \succ q_1$. The corresponding reduced
Gr\"{o}bner  basis  of $I$ is $F = \{ q_1^3 q_4 - q_2^2 q_3,\, q -
q_1,\,  q_1^2  x_2  -  q_2,\,  q_2  q_3  x_2 - q_1 q_4 ,\, q_1 q_3
x_2^2-q_4,\,  q_1  q_2  x_3-q_3,\, q_2^2 x_3- q_1 q_3 x_2,\, q_1^2
q_4  x_3-  q_3^2 x_2,\, q_2 q_4 x_3- q_3^2 x_2^2,\, q_1 q_4^2 x_3-
q_3^3  x_2^3,\,  q_4^3  x_3-  q_3^4  x_2^5  \}$.
By the  elimination
property  of  Gr\"{o}bner  bases $G := F \cap k[q_1,
q_2,  q_3,  q_4]  =  \{  q_1^3  q_4  -  q_2^2  q_3 \}$ is
a Gr\"{o}bner  basis  for  the  elimination ideal $J = I \cap k[q_1,
q_2,  q_3, q_4]$.
Collecting  leading  terms  of  $G$
gives  $L  =  \{  q_1^3  q_4  \}$.
Since $G$ is a Gr\"obner basis of $J$ we know that $I_\preceq(J) =
\langle  q_1^3  q_4  \rangle$.  The  Hilbert-Poincar\'{e} series
of  $k[q_1, q_2, q_3, q_4] / \langle q_1^3 q_4 \rangle$ gives
\begin{equation}\label{eq:4}
\mathcal{M}_{[4]}(q) = \frac{1 - q^7}{(1-q) (1-q^2) (1-q^3)
(1-q^4)}.\end{equation}
It is clear that we may replace
any occurrence of the partition $4+1+1+1$ in
a multinomial coefficient
by $3+2+2$ without changing
the value of the multinomial coefficient.
Therefore, there are at most
as many multinomial coefficients
as there are partitions
avoiding $4 + 1 + 1 + 1$.
Equation~(\ref{eq:4}) states that
this upper bound gives in fact the exact number in the case
of $S=\{1,2,3,4\}$.

Note that all denominators in the partial fraction
decomposition
\begin{multline}\label{eq:M4pfd}
\mathcal{M}_{[4]}(q) =
-\frac{7}{24} \frac{1}{(q-1)^3}
- \frac{77}{288} \frac{1}{(q-1)}
+ \frac{1}{16} \frac{1}{(q+1)^2}  + \\
+ \frac{1}{32} \frac{1}{(q+1)}
+ \frac{1}{9} \frac{(q+2)}{(q^2 + q + 1)}
+ \frac{1}{8} \frac{(q + 1)}{(q^2 + 1)}
\end{multline}
of (\ref{eq:4}) are powers of cyclotomic polynomials
$C_j(q)$ where $j$ divides an element of $S = \{1,2,3,4\}$.
We rewrite this as
\begin{multline}\label{eq:M4pfderweitert}
\mathcal{M}_{[4]}(q) =
\frac{7}{24} \frac{1}{(1 - q)^3}
+ \frac{77}{288} \frac{1}{(1 - q)}
+ \frac{1}{16} \frac{(1-q)^2}{( 1 - q^2 )^2}  + \\
+ \frac{1}{32} \frac{(1-q)}{(1 - q^2)}
+ \frac{1}{9} \frac{(2 - q - q^2)}{(1 - q^3)}
+ \frac{1}{8} \frac{(1 + q - q^2 - q^3)}{(1 - q^4)}.
\end{multline}
in order to use the binomial series
$(1-z)^{-a-1} = \sum_{n=0}^\infty \binom{a+n}{a} z^n$.
The result is
\begin{multline}\label{eq:M4n}
M_{[4]}(n) =
\frac{7}{48} n^2
+ \left(\frac{1}{16} [1,-1](n) +  \frac{7}{16} \right) n + \\
+       \frac{1}{8} [1,1,-1,-1](n)
      + \frac{1}{9} [2,-1,-1](n)
      + \frac{3}{32}[1,-1](n)
      + \frac{161}{288}
\end{multline}
where $[a_0, a_1, \dots, a_m](n) := a_j$ for $n \equiv j (m)$.
Similar  computations  show that
\begin{equation}\label{eq:5}
\mathcal{M}_{[5]}(q) =
\frac{1 - q^7}{(1-q) (1-q^2) \dots (1-q^5)},\end{equation}
\begin{equation}\label{eq:6}
\mathcal{M}_{[6]}(q) = \frac{1 - q^7 - q^8 - q^{10} +
q^{12} + q^{13}} {(1-q) (1-q^2) \dots
(1-q^6)},\end{equation} and
\begin{equation}\label{eq:7}
\mathcal{M}_{[7]}(q) = \frac{1 - q^7 - q^8 - q^{10} + q^{12} +
q^{13}} {(1-q) (1-q^2) \dots
(1-q^7)}.\end{equation}
It is no coincidence that the numerators of (\ref{eq:6}) and (\ref{eq:7})
agree (Theorem~\ref{thm:SuP'}).

\section{Upper Bounds}
\label{sec:uppers}  Trivially,  $M(n) \leq p(n)$. Our goal is to
find sharper upper bounds.
\begin{lem}\label{lem:upper}
Assume  $S'\subseteq  S$.

Let $\tilde{I}$ be the ideal of $k[S,q,x]$ generated by the set of
polynomials $\{ q_j - \h(q_j) \,:\, j\in S'\}$. Let $\tilde{J}$ be
the ideal generated by $\tilde{I} \cap k[S']$ in the ring $k[S]$.
Let $U_{S,\,S'}(n) := H_{k[S] / \tilde{J}}(n).$ Then
\begin{enumerate}
\item
$M_S(n) \leq U_{S,\,S'}(n)$.
\item
We have
\begin{equation}\label{eq:h1}
\sum_n U_{S,\,S'}(n) q^n
= \frac{f_{S'}(q)}{\prod_{j\in S} (1-q^j)}
\end{equation}
where $f_{S'}(q)$ is defined by
\begin{equation}\label{eq:h2}
\sum_n M_{S'}(n) q^n
= \frac{f_{S'}(q)}{\prod_{j\in S'} (1-q^j)}.
\end{equation}
\end{enumerate}
\end{lem}
\begin{pf}
We prove the first statement. Let $I$ be the ideal of $k[S,q,x]$
generated by the set of polynomials $\{ q_j - \h(q_j) \,:\, j\in
S\}$ and let $J = I \cap k[S]$. Since $\tilde{J}$ is a $k$-vector
subspace of $J$ we have
\begin{equation}
\dim_k k[S]_n \cap J \geq \dim_k k[S]_n \cap \tilde{J} \end{equation}
and therefore
\begin{equation} \dim_k (k[S] / J)_n
\leq \dim_k (k[S] /\tilde{J})_n\end{equation}
i.e.
\begin{equation} M_S(n) \leq U_{S,\,S'}(n).\end{equation}

To prove the second statement, let $I'$ be the ideal generated by
$\{q_j - \h(q_j) \,:\, j\in S'\}$ in the ring $k[S',q,x]$ and let
$J' := I' \cap k[S']$. Since the ideals $\tilde{J}$ and $J'$ are
generated by the same set of polynomials (albeit in different
rings), the Hilbert functions $U_{S,S'}(n) = H_{k[S] /
\tilde{J}}(n)$ and $M_S'(n) = H_{k[S']/ J'}(n)$ correspond in the
way claimed by (\ref{eq:h1}) and (\ref{eq:h2}). \qed\end{pf} To
get upper bounds for $M(n)$, we use the preceding Lemma in the
special case $S=\mathbb{N}$ getting:
\begin{thm}\label{thm:Mleq}\label{thm:Mpp7}
For any $S'$ we have
\begin{equation}
M(n) \leq [q^n] \frac{f_{S'}(q)}{\prod_{j=1}^{\infty} (1-q^j)}
\end{equation}
(where $[q^n]\mathcal{A}(q)$ denotes the the coefficient of $q^n$ in the
power series expansion of $\mathcal{A}(q)$).
For instance, the cases $S'=[4]$ and $S'=[6]$ yield the bounds
\begin{equation}
        M(n)        \leq        p(n)       -       p(n-7),
\end{equation}
and
\begin{equation}
  M(n) \leq p(n) - p(n-7) - p(n-8) - p(n-10) + p(n-12) + p(n-13).
\end{equation}
\end{thm}
Note that a  direct proof of $M(n) \leq p(n) - p(n-7)$ could
be given by exploiting the equivalence of the partitions $4+1+1+1$
and $3+2+2$ in the sense of Equation~(\ref{eq:3224111}).

The bound $M(n) \leq p(n) - p(n-7)$ is good enough to imply:
\begin{thm}\label{thm:Mp0}
$M(n)=o(p(n))$,  i.e.  $\lim_{n\rightarrow\infty}  M(n)/p(n) = 0$.
\end{thm}

\begin{pf}
Due  to  the  monotonicity  of  $p(n)$  and the fact that the unit
circle is the natural boundary for
\begin{equation}
\sum_{n=0}^{\infty} p(n) q^n
        = \prod_{n=1}^{\infty} \frac{1}{1 - q^n}\;, \end{equation}
we see that
\begin{equation}
\lim_{n\rightarrow   \infty}   \frac{p(n   -   7)}{p(n)}   =  1\,.
\end{equation} Hence
\begin{equation}
        0  \leq  \lim_{n\rightarrow \infty} \frac{M(n)}{p(n)} \leq
        \lim_{n\rightarrow  \infty}  \frac{p(n)  - p(n - 7)}{p(n)}
        = 1 - 1 = 0\,,
\end{equation}
which proves Theorem \ref{thm:Mp0}.
\qed\end{pf}

\section{Lower Bounds} \label{sec:lowers}
Recall that $M(n) \geq p_{\Primes}(n)$
(Theorem~\ref{lower-bound-thm}).
The numbers given on page \pageref{numbers} suggest
that $M(n)$ grows much faster than $p_{\Primes}(n)$.
We will prove that this is indeed the case:
$\lim_{n\rightarrow\infty}  p_{\Primes}(n)/M(n) = 0$
(Theorem \ref{thm:pM0})
and we will give better lower bounds for $M(n)$.

Let us write $S < P$ if each element of $S$ is less than each element
of $P$.
We need the following generalization of
Lemma~\ref{NeverEquivalent}:
\begin{lem}\label{lem:pp}
Assume $S < P$ where $P$ is a set of primes. Let $s$ and $s'$ be
any two power products in $k[S]$ and let $p$ and $p'$ be distinct
power products in $k[P]$. Then $\h(s p) \neq \h(s' p')$.
\end{lem}

In the case $S=\emptyset$, Lemma~\ref{lem:pp} states that distinct
partitions $p$ and $p'$ into primes yield different multinomial
coefficients: $\h(p) \neq \h(p')$. Lemma~\ref{lem:pp} can be
proved by the same induction argument as
Lemma~\ref{NeverEquivalent}.

\begin{lem}\label{lem:SuP}
Assume $S < P$ where $P$ is a set of primes. Define $\h$ on $k[S
\cup P]$ by (\ref{eq:phidef}). Then $\ker \h$ is generated, as an
ideal of $k[S \cup P]$, by $\ker \h \, \cap \, k[S]$.
\end{lem}
\begin{pf}
Let $f \in \ker \h$. Since $k[S \cup P] = k[S][P]$, we can $f$ as
a finite sum $f = \sum_s \sum_p c_{s,p}\, s p$ indexed by power
products $s$ and $p$ from $k[S]$ and $k[P]$ respectively, with
coefficients $c_{s,p} \in k$. As $f \in \ker \h$, $\sum_s \sum_p
c_{s,p}\, \h(s p) = 0.$ By Lemma~\ref{lem:pp}, this implies
$\sum_s c_{s,p}\, \h(s p) = 0 $ for arbitrary but fixed $p$.
Cancelling $\h(p)$ from this equation shows that $\h(f_p)=0$ where
$f_p := \sum_s c_{s,p}\, s$. In this way we succeed in writing $f$
as $f = \sum_p f_p p$ where each $f_p$ is in $\ker \h \cap k[S]$.
\qed\end{pf} As an immediate consequence of Lemma~\ref{lem:SuP} we
get:
\begin{thm}\label{thm:SuP'}
Assume $S < P$ where $P$ is a set of primes.
Then
\begin{equation}
\mathcal{M}_{S \cup P}(q) = \mathcal{M}_{S}(q) / \prod_{j\in P} (1-q^j).
\end{equation}
\end{thm}
As a first application of Theorem~\ref{thm:SuP'}, we
count multinomial coefficients with lower entries which
are either prime or equal to $1$:
\begin{equation}\label{eq:M1P}
\mathcal{M}_{\{1\} \cup \Primes}(q) =
\frac{1}{(1-q) \prod_{j \in \Primes}(1-q^j) },
\end{equation}
which
allows for improving
Theorem~\ref{lower-bound-thm}:
\begin{thm}\label{thm:pM0}
We have
\begin{equation}
\lim_{n\rightarrow\infty}
p_{\Primes}(n) / M_{\{1\} \cup \Primes}(n) = 0
\end{equation}
and therefore
$\lim_{n\rightarrow\infty}  p_{\Primes}(n)/M(n) = 0.$
\end{thm}
\begin{pf}
Let $A(n) := M_{\{1\} \cup \Primes}(n)$.
Due  to  the  monotonicity  of  $A(n)$
and the fact that the unit
circle is the natural boundary for we see that
\begin{equation}\lim_{n\rightarrow \infty} A(n-1) / A(n)
=  1.\end{equation}
By (\ref{eq:M1P}),
\begin{equation}
p_{\Primes} (n) =
A(n) - A(n-1).
\end{equation}
Therefore,
\begin{equation}
        0  \leq  \lim_{n\rightarrow \infty}
         \frac{p_\Primes(n)}{A(n)} \leq
        \lim_{n\rightarrow  \infty}  \frac{A(n) - A(n-1)}{A(n)}
        = 1 - 1 = 0\,,
\end{equation}
which proves Theorem \ref{thm:pM0}.
\qed\end{pf}

Let $L_S(n) := M_{S \cup \Primes}(n)$; clearly, $L_S(n)$ is a
lower bound for $M(n)$.
Theorem~\ref{thm:SuP'} allows us deduce
\begin{equation}\label{eq:L45}
\mathcal{L}_{[4]}(q) = \mathcal{L}_{[5]}(q) = \frac{1 - q^7}
{
\prod_{j \in [4] \cup \Primes} (1-q^j)
}
\end{equation}
and
\begin{equation}\label{eq:L67}
\mathcal{L}_{[6]}(q) = \mathcal{L}_{[7]}(q) =
\frac{1 - q^7 - q^8 - q^{10} +
q^{12} + q^{13}}
{
\prod_{j \in [6] \cup \Primes} (1-q^j)
}\end{equation}
from the Equations (\ref{eq:4}) -- (\ref{eq:7}); some values
of $L_{[4]}(n)$  are listed on page \pageref{numbers}.

\section{The Irreducible Pairs} \label{sec:in}

An  {\em  irreducible  pair}  is  a pair of partitions of $n$ that
yield  the  same  multinomial  coefficient  but  have  no parts in
common. For example,
\begin{equation}\label{eq:pair4111}
    (4,1,1,1)
        \text{ and }
    (3,2,2)
\end{equation}
is an irreducible pair.

It  turns  out that there are infinitely many irreducible
pairs of partitions.  The following is a partial list:
Generalizing (\ref{eq:pair4111}) we see that
\begin{equation}
        (2^m,\underbrace{1,1,\dots,1}_{2m-1}\,) \text{ and }
        (2^m - 1, \underbrace{2,2,\dots,2}_{m}\,)
\end{equation} form an irreducible pair of partitions of $2^m + 2m - 1$.
More generally, for any integers $a \geq 2$ and $m \geq 1$
the partitions
\begin{equation}
\label{eq:pn}
        (a^m, \underbrace{a-1,a-1,\dots,a-1}_{m},
        \underbrace{1,1,\dots,1}_{m-1}\,) \text{ and }
        (a^m - 1, \underbrace{a,a,\dots,a}_m \,)
\end{equation} form an irreducible pair of partitions of $a^m + am - 1$.

The pair
\begin{equation} \label{eq:pair611}
    (6,1,1) \text{ and } (5,3)
\end{equation}
can be generalized to irreducible pairs
\begin{equation}\label{eq:jfactorial}
    (j!,\, \underbrace{1,\dots,1}_{(j-1)}\,)
        \text{ and }
    (j!-1,\,j)
\end{equation}
of partitions of $(j! + j - 1)$ for $j \geq 3$.

From any two irreducible pairs we can get a third one
by combining them in a natural way.
For instance, combining $a$ copies of (\ref{eq:pair611})
with $b$ copies of (\ref{eq:pair4111}) gives
the pair (\ref{eq:ab}) which is used in the proof below.

The above examples show that $i(n)$ is positive infinitely often.
Indeed we have:
\begin{thm} \label{thm:n56}
$i(n) \geq \frac{n}{56} - 1$.
\end{thm}
\begin{pf}  For each pair of non-negative integers $a$ and $b$
satisfying \begin{equation} \label{eq:37}
        8a + 7b = n\,,
\end{equation} we see that \begin{equation} \label{eq:ab}
        (\underbrace{6,\dots,6}_a,\underbrace{4,\dots,4}_b,
        \underbrace{1,\dots,1}_{2a+3b}) \text{ and }
        (\underbrace{5,\dots,5}_a,\underbrace{3,\dots,3}_{a+b},
        \underbrace{2,\dots,2}_{2b})
\end{equation} forms a new irreducible pair of partitions of $n$. Consequently
$i(n)$ is at least as large as the number of non-negative
solutions of the linear Diophantine equation (\ref{eq:37}).

Now the segment of the line $8a + 7b = n$ in the first quadrant is
of length $n\sqrt{113}\big/56$.  Furthermore from the full
solution of the linear Diophantine equation we note that the
integral solutions of (3.7) are points on this spaced a distance
$\sqrt{113}$ apart.  Hence in the first quadrant there must be at
least \begin{equation}
        \left\lfloor \frac{n\sqrt{113}/56}{\sqrt{113}}\right\rfloor
        = \left\lfloor \frac{n}{56} \right\rfloor > \frac{n}{56} - 1
\end{equation} such points.  Therefore \begin{equation}
        i(n) > \frac{n}{56} - 1\,.
\end{equation}
\qed\end{pf}

Theorem~\ref{thm:n56} shows that $i(n) > 0$ for all $n \geq 56$.
Direct computation shows that $i(n) > 0$ for all $n > 7$ with the
exception of $n = 9,11$ and $12$.

\section{Further Problems}
Clearly we have only scratched the surface concerning the order of
magnitude of $M_k(n)$, $M(n)$ and $i(n)$.  We have computed tables
of the functions, and based on that evidence we make the following
conjectures.
\begin{conj}\label{con:p*}
$M(n)  \geq  p^*  (n)$ for $n \geq 0$, where $p^*(n)$ is the total
number of partitions of $n$ into parts that are either $\leq 6$ or
multiples of $3$ or both.
\end{conj}
\begin{conj}\label{con:C}
There     exists    a    positive    constant    $C$    so    that
\begin{equation}
        \lim_{n\rightarrow\infty}  \frac{\log M(n)}{\sqrt{n}} = C.
\end{equation}
\end{conj}
If  $C$  exists  and  if  Conjecture  \ref{con:p*}  is  true, then
\cite[Th. 6.2, p.89]{A76}
\begin{equation}
        \frac{\pi}{3}  \sqrt{2}  \leq C \leq \pi \sqrt{\frac23}\;.
\end{equation}
\begin{conj}\label{con:Ck}
Let $C_k$ be the infimum of the quotients $M_k(n) / p_k(n)$ where $n$ ranges
over the natural numbers.
Then $C_k > 0$ for all natural numbers $k$. Moreover,
$C_k$ is a strictly decreasing function of $k$ for $k
\geq 3$ and $C_k \to 0$ as $k \to \infty$.
\end{conj}
\begin{conj}\label{con:phash}$M(n) \leq p^{\#}(n)$ for $n \geq 0$ where
$p^{\#}(n)$ is the total number of partitions of $n$ into parts
that are either $\leq 7$ or multiples of $3$ or both.
\end{conj}
Conjecture \ref{con:phash} together with Conjecture \ref{con:p*}
allows us to replace Conjecture \ref{con:C} with
\begin{conj}
\begin{equation}
        \lim_{n\rightarrow \infty} \frac{\log M(n)}{\sqrt{n}} =
        \frac{\pi}{3} \sqrt{2}\,.
\end{equation}
\end{conj}

{\em Acknowledgements:} {\small We thank Anna Bigatti for expert
advice on \cocoa \cite{CocoaSystem} and Bogdan Matasaru for
rewriting a program for computing $M(n)$ in the programming language
{\tt C}.}

\begin{table}
\begin{tabular} {rrrrrrrr}
$n$ & $p_{\Primes}(n)$ & $L_{[4]}(n)$ & $p^*(n)$ & $M(n)$
&
$p^{\#}(n)$ & $U_{\mathbb{N}^+,[4]}(n)$ & $p(n)$\\
\hline
$0$ & $1$ & $1$ & $1$ & $1$ & $1$ & $1$ & $1$ \\
$10$ & $5$ & $30$ & $36$ & $36$ & $39$ & $39$ & $42$ \\
$20$ & $26$ & $232$ & $357$ & $366$ & $445$ & $526$ & $627$ \\
$30$ & $98$ & $1102$ & $2064$ & $2131$ &
  $2875$ & $4349$ & $5604$ \\
$40$ & $302$ & $4020$ & $8853$ & $9292$ &
  $13549$ & $27195$ & $37338$ \\
$50$ & $819$ & $12405$ & $31639$ & $33799$ &
  $52321$ & $140965$ & $204226$ \\
$60$ & $2018$ & $34016$ & $99245$ & $107726$ &
  $175426$ & $636536$ & $966467$ \\
$70$ & $4624$ & $85333$ & $281307$ & $310226$ &
  $527909$ & $2582469$ & $4087968$ \\
\hline
\end{tabular}
\label{numbers}
\end{table}

\nocite{A94}
\bibliographystyle{elsart-num}
\bibliography{mn}
\end{document}